\documentstyle[12pt,twoside]{amsart}
\title{Minimal model theorem for toric divisors}
\address{Department of Mathematics\\ Tokyo Institute of Technology\\
Oh-Okayama, Meguro, Tokyo 152, JAPAN}
\author{Shihoko Ishii}

\newcommand{\bZ}{{\Bbb Z}}
\newcommand{\bR}{{\Bbb R}}
\newcommand{\bN}{{\Bbb N}}
\newcommand{\bC}{{\Bbb C}}
\newcommand{\bQ}{{\Bbb Q}}
\newcommand{\bP}{{\Bbb P}}

\newcommand{\D}{{\Delta}}

\newcommand{\bp}{{\bold p}}
\newcommand{\bq}{{\bold q}}

\newcommand{\bn}{{\bold n}}
\newcommand{\bm}{{\bold m}}

\newcommand{\NR}{{N_{\bR}}}
\newcommand{\MR}{{M_{\bR}}}
\newcommand{\boxh}{{\Box_h}}
\newcommand{\boxk}{{\Box_k}}
\newcommand{\tnd}{{T_N(\D)}}
\newcommand{\tns}{{T_N(\Sigma)}}
\newcommand{\xs}{{X(\Sigma)}}
\newcommand{\epm}{{\epsilon^{(m)}}}
\newcommand{\km}{{k^{(m)}}}
\newcommand{\tsig}{{\tilde \Sigma}}
\newtheorem{thm}{Theorem}[section]

\newtheorem{lem}[thm]{Lemma}
\newtheorem{cor}[thm]{Corollary}
\newtheorem{prop}[thm]{Proposition}
\newtheorem{claim}[thm]{Claim}

\theoremstyle{definition}
\newtheorem{defn}[thm]{Definition}

\newtheorem{say}[thm]{}
\newtheorem{exmp}[thm]{Example}

\theoremstyle{remark}
\newtheorem{case}{Case}

\begin{document}
\bibliographystyle{amsplain}
\maketitle
\pagestyle{myheadings}
\markboth{\hfill SHIHOKO ISHII\hfill}{\hfill MINIMAL MODEL THEOREM \hfill}

\begin{abstract}
  Minimal model conjecture for a proper variety $X$ is that 
  if $\kappa(X)\geq 0$, then $X$  has a minimal model with the abundance
  and if $\kappa =-\infty$, then $X$  is birationally equivalent
  to a variety
  $Y$  which has a fibration $Y \to Z$ with $-K_Y$ relatively 
  ample.
  In this paper, we prove this conjecture for a $\D$-regular
  divisor on a proper toric variety
  by means of successive contractions of extremal rays and flips of
  ambient toric variety.
  Furthermore, for such a divisor $X$ with $\kappa(X)\geq 0$
  we construct a projective minimal model with the abundance 
  in a different way;
  by means of "puffing up" of the polytope,
  which gives an algorithm of a construction of a minimal model.
\end{abstract}

\setcounter{section}{-1}

\section{\bf Introduction}

  Let  $k$  be an algebraically closed field of arbitrary 
  characteristic. 
  Varieties in this paper are all defined over  $k$.
  Let  $X$  be a proper algebraic variety.
   A proper algebraic variety  $Y$  is called a minimal model of  $X$,
  if (1) $Y$  is birationally equivalent to $X$,
   (2)  $Y$  has at worst terminal singularities and
   (3)  the canonical divisor  $K_Y$  is nef.
  A minimal model  $Y$  is said to have the abundance if the linear
  system  $|mK_Y|$  is basepoint free for sufficiently large $m$.
  The minimal model  conjecture states:
  an  arbitrary proper variety with $\kappa\geq 0$ has a minimal model
  with the abundance 
  and an arbitrary proper variety with $\kappa=-\infty$ has a 
  birationally
  equivalent model  $Y$  with at worst terminal singularities 
  and a fibration $Y\to Z$ to a lower dimensional variety
  with $-K_Y$ relatively ample.

  The conjecture holds true for 2-dimensional case which is known 
  as a classical result.
  For 3-dimensional case the conjecture for $k=\bC$
  is proved by Mori \cite{M}
  and Kawamata \cite{kawamata}, while it is not yet proved for higher
  dimensional case.
  As a special case of higher dimension,
  Batyrev \cite{Batyrev} proved, among other results,
  the existence of a minimal model
  for a $\D$-regular anti-canonical divisor  of a Gorenstein Fano 
  toric variety
  $\tnd$.
  In this paper
  we prove the minimal model conjecture for every $\D$-regular
  divisor $X$ on a toric variety of arbitrary dimension
  by means of successive contractions of extremal rays
  and flips.
  Furthermore for such a divisor with $\kappa \geq 0$,
  we construct a projective minimal model with the abundance
  in a different way;
  by means of "puffing up" of the polytope corresponding to the 
  adjoint divisor.
  By this method one can concretely construct a projective minimal
  model.
  As a corollary,
  for a field  $k$  of characteristic 0,
  the minimal model conjecture holds for  
  a general member 
  of a basepoint free linear system on a proper  
  toric variety over $k$. 
  The half of this work was done during the author's stay at the 
  Johns
  Hopkins University on April 1996.
  She expresses her gratitude to Professors Shokurov and Kawamata
  who made her stay possible.
  She is also grateful to the Johns Hopkins University for their 
  hospitality.
  She would like to thank Professor Reid who gave  useful suggestions
  and 
  Professor Batyrev who called her attention
  to this problem and pointed out an error of the first draft of this
  paper.

\section{\bf The minimal model theorem for toric divisors}

\begin{defn}
(\cite{Batyrev}) 
  A  divisor  $X$  of a toric variety $\tnd$ defined by a fan
  $\D$ is called
  $\D$-regular,
  if for every  $\tau \in \D$ the  intersection  
  $X\cap orb( \tau)$ is either a smooth
  divisor of  $orb( \tau)$  or empty.
\end{defn}

\begin{defn}
  Let $V$  and $V'$  are toric varieties defined by fans $\D$ and
  $\D'$ respectively  
  and  $f: V' - \to V$   a toric birational map: i.e. $\D'$  is 
  obtained by successive subdivisions and converse of subdivisions
  from $\D$.
  Let $T$  be the maximal orbit  in $V$.
  If an irreducible divisor  $X$ on $V$  satisfies 
  $X\cap T\neq \phi$,
  the divisor  $X'=\overline{f^{-1}(X\cap T)}$ on $V'$
  is called the proper transform of  $X$ on $V'$.
\end{defn}

\begin{defn}
  Let  $X$ a divisor on a normal variety  $V$  such that  $K_V+X$
  is a $\bQ$-Cartier divisor and $f:V'\to V$ a birational morphism.
  Let  $X'$ be the proper transform of $X$.
  If 
$$K_{V'}+X'=f^*(K_V+X)+\sum_i a_iE_i,$$
  where $E_i$'s are the exceptional divisors of $f$,
  then  $a_i$  is called the discrepancy of $K_V+X$  at  $E_i$
\end{defn}

\begin{defn}
   Let  $V$  be a toric variety defined by a simplicial fan $\D$
  and $X$ an irreducible divisor on  $V$.
  The divisor $K_V+X$  is called terminal, if the following hold:

(1) there exists a morphism $f:V'=T_N(\D')\to V$ corresponding to
  a non-singular subdivision  $\D'$ of $\D$ ($\D'\neq \D$) 
  such that the proper transform $X'$  of $X$  on $V'$  is 
  $\D'$-regular, 
  in particular $X\cap T\neq \phi$ for the maximal orbit $T$
  in $V$,
  and 

(2) for every such morphism  as in (1) the discrepancy of $K_V+X$ 
  at every exceptional divisor on $V'$  is positive.
\end{defn}

\begin{lem}
\label{terminal}
  If $V=\tnd $ is non-singular and an irreducuble divisor $X$ 
  on $V$  is 
  $\D$-regular,
  then  $K_V+X$  is terminal
\end{lem}

\begin{pf}
  For every non-singular subdivision  $\D'$  of  $\D$, where 
  $\D'\neq \D$,
  the proper transform $X'$ of $X$ by the corresponding morphism
  $f:V'=T_N(\D')\to V$ is $\D'$-regular by 3.2.1 of \cite{Batyrev}.
  Since $X'=f^*X$  and $K_{V'}=f^*K_V+\sum_ia_iE_i$, where $a_i>0$
  for every exceptional divisor $E_i$  on $V'$,
  it follows that the discrepancy of $K_{V}+X$  at each $E_i$ is 
  positive.
\end{pf}

\begin{prop}
   Let  $V$  be a toric variety defined by a simplicial fan $\D$
  and $X$ an irreducible  divisor on  $V$.
  Then the divisor $K_V+X$  is  terminal if and only if the following 
  hold:

(i) there exists a morphism $f:V'=T_N(\D')\to V$ corresponding to
  a non-singular subdivision  $\D'$ of $\D$ ($\D'\neq \D$) 
  such that the proper transform $X'$  of $X$  on $V'$  is 
  $\D'$-regular. 

(ii) for one such morphism as in (i) the discrepancy of $K_V+X$ 
  at every exceptional divisor on $V'$  is positive.
\end{prop}

\begin{pf}
  Let $f:V'=T_N(\D')\to V$ be the morphism satisfying the condition
  (i)  and (ii) and $g:V''\to V$  be another morphism satisfying (i).
  Take a nonsingular toric variety $\tilde V$ which dominates both
  $V'$ and $V''$.
  Then by \ref{terminal}, $K_{V'}+X'$ is terminal.
  Therefore  the discrepancy of $K_V+X$ at every exceptional divisor
  on $\tilde V$  is positive which yields the positivity of it
  at every exceptional divisor on $V''$.
\end{pf}

\begin{lem}
\label{terminal sing}
  Let  $V$  be a toric variety defined by a simplicial fan $\D$
  and $X$ an irreducible divisor on  $V$.
  If the divisor $K_V+X$  is  terminal,
  then $V$ has at worst terminal singularities.
\end{lem}
\begin{pf}
This follows from the fact that a discrepancy of $K_V$  is 
  greater than or equal to that of $K_V+X$. 
\end{pf}

  Here we summerize the results of Reid (\cite{Reid}) which are used
  in this section.

\begin{prop}
\label{reid}
(\cite{Reid})
  Let $V$  be the toric variety defined by a proper simplicial fan 
  $\D$.
  
(i) $NE(V)=\sum_{i=1}^r\bR_{\geq 0}[\ell_i]$,
   where $\ell_i$'s are 1-dimensional strata on $V$.
  Here each $\bR_{\geq 0}[\ell_i]$  is called an extremal ray.

(ii) For every extremal ray  $R$ there exist a toric morphism
  $\varphi_R:V\to V'$ which is an elementary contraction in the sense
  of Mori theory: $\varphi_R{\cal O}_V={\cal O}_{V'}$
  and $\varphi _RC=pt$ if and ony if $[C]\in R$.
  Let $A\subset V$ and $B\subset V'$ be the loci on which
  $\varphi_R$  is not an isomorphism, then $\varphi_R|_A:A \to B$
  is a flat morphism and all of whose fibers are weighted projective
  spaces of the common dimension.  

(iii)  If $\varphi_R:V\to V'=T_N(\D')$  is birational and not 
  isomorphic in codimension one,
  then the exceptional set of $\varphi_R$  is an irreducible 
  divisor and $\D'$  is proper simplicial.
  Here this $\varphi_R$  is called a  divisorial contraction.

(iv)  If $\varphi_R:V\to V'=T_N(\D')$  is  
  isomorphic in codimension one, then there exists the following
  commutative diagram:
 $$\matrix
  &              & \tilde V &                & \\
  &\psi \swarrow &          &\searrow \psi_1 & \\
V &              &          &                &V_1=T_N(\D_1)\\
  &\varphi_R \searrow&      &\swarrow \varphi_1&  \\
  &              & V'       &                 &   
\endmatrix$$
  such that $\D_1$  is proper simplicial, 
  $\D_1(1)=\D(1)$,
  all morphisms are elementary 
  contractions of extremal rays,
  $\psi$ and $\psi_1$ are birational morphisms with the exceptional
  divisor $D$,
  $\varphi_R$ and $\varphi_1$ are birational morphisms with the 
  exceptional sets
  $\psi(D)$ and $\psi_1(D)$ respectively,
  and identifying $N_1(V)$ and $N_1(V_1)$, $-R$  is an extremal ray in 
  $NE(V_1)$  and $\varphi_1=\varphi_{-R}$.
  Here the birational map $\varphi_1^{-1}\circ \varphi_R:V- \to V_1$
  is called a  flip.
\end{prop}

\begin{lem}
\label{elementary}
  Let  $V$  be a toric variety defined by a proper simplicial fan
  $\D$ and $X$ an irreducible divisor such that $K_V+X$  is terminal.
  Let  $R$  be an extremal ray such that $(K_V+X)R<0$.
  Then the following hold:

(i) if $\varphi_R:V\to V'=T_N(\D')$ is a divisorial contraction,
  then $K_{V'}+X'$  is terminal, where $X'$  is the proper transform
  of  $X$ on $V'$;

(ii) let $\varphi_R:V\to V'$ be isomorphic in codimension one;
  for the diagram 
 $$\matrix
  &              & \tilde V &                & \\
  &\psi \swarrow &          &\searrow \psi_1 & \\
V &              &          &                &V_1=T_N(\D_1)\\
  &\varphi_R \searrow&      &\swarrow \varphi_1&  \\
  &              & V'       &                 &   
\endmatrix$$
of (iv), \ref{reid},
  let $X_1$  be the proper transform of $X$  on $V_1$,
  $D$ the exceptional divisor of $\psi$ and $\psi_1$,
  $\alpha$ the discrepancy of $K_V+X$  at $D$
  and $\alpha'$ the discrepancy of $K_{V_1}+X_1$  at $D$;
  then $\alpha <\alpha'$  and $K_{V_1}+X_1$  is terminal.
\end{lem}

\begin{pf}
  For the proof of (i), first one should remark that $V'$  is
  $\bQ$-factorial, because $\D'$  is simplicial.
  Let  $E$  be the exceptional divisor for $\varphi_R$.
\begin{claim}
$ER<0$.
\end{claim}
  For the proof of the claim, take an irreducible divisor $H$ on $V'$
  such that $H\supset \varphi_R(E)$.
  Then  $\varphi^*H=[H]+aE$ with $a>0$, 
  where  $[H]$  is the proper transform of  $H$ on $V$.
  Since $ (\varphi_R^*H) R=0$ and $[H]R>0$, it follows that 
  $aER<0$ which completes the proof of the claim.

  Denote $K_V+X$  by $\varphi_R^*(K_{V'}+X')+bE$,
  then $b>0$.
  In fact, by $(K_V+X)R<0$, $\varphi_R^*(K_{V'}+X')R=0$ and
   $ER<0$, it follows that $b>0$.
  Let  $\overline \D$  be a non-singular subdivision of $\D$
  such that the proper transform $\overline X$  of  $X$  on 
  $\overline V=T_N(\overline \D)$ is $\overline \D$-regular.
  Since  $K_V+X$   is terminal, the discrepancy of $K_V+X$ at 
  every exceptional divisor for $\overline V \to V$  is positive.
  By this, and $b>0$, it follows that the discrepancy of $K_{V'}+X'$
  at every exceptional divisor for $\overline V \to V'$ is positive.
  For the proof of (ii),
  take a curve  $\ell$  on $\tilde V$ such that $\psi_1(\ell)=pt$
  and $\psi(\ell)\neq pt$.
  This is possible,
  because if a curve contracted by both $\psi$ and $\psi_1$ exists,
  then the extremal rays corresponding to $\psi$ and $\psi_1$
  coincide which implies $V\simeq V_1$ and $\varphi_R=\varphi_1$
  a contradiction to $\varphi_1=\varphi_{-R}$  in (iv) of \ref{reid}.
  For this $\ell$, one can prove that $D\ell <0$  in the same way
  as in the claim above.
  Now as $\psi_*(\ell)$  is contracted to a point by $\varphi_R$,
  $[\psi_*(\ell)]\in R$, therefore
  $\psi^*(K_V+X)\ell=(K_V+X)\psi_*(\ell)<0$.
  By intersecting  $\ell$ with $K_{\tilde V}+\tilde X=
  \psi^*(K_V+X)+\alpha D$,
  we obtain 
  $$(K_{\tilde V}+\tilde X)\ell< \alpha D\ell.$$
  Here the left hand side is $\psi_1^*(K_{V_1}+X_1)\ell + \alpha'D\ell$,
  and $\psi_1^*(K_{V_1}+X_1)\ell=0$ because of the definition of 
  $\ell$.
  This proves  that $\alpha<\alpha'$.
 To prove the last statement,
  take a non-singular subdivision $\tilde {\tilde \D}$ of $\tilde \D$
  such that the proper transform $\tilde {\tilde X}$  of  $\tilde X$
  is  $\tilde {\tilde \D}$-regular.
  Let $\lambda:\tilde {\tilde V}=T_N(\tilde {\tilde \D})\to \tilde V$
  be the corresponding morphism.
  Then $K_{\tilde {\tilde V}}+\tilde {\tilde X}=
   \lambda^*\psi^*(K_V+X)+\sum_i\beta_iE_i$,
  where $\beta_i>0$ for every exceptional divisor $E_i$, 
  because $K_V+X$  is terminal.  
  Now by substituting $\psi^*(K_V+X)=\psi_1^*(K_{V_1}+X_1)+
  (\alpha'-\alpha)D$
  into the equality above,
  the discrepancy of $K_{V_1}+X_1$  at every exceptional divisor
  on $\tilde {\tilde V}$ turns out to be positive.
\end{pf}

\begin{thm}
\label{theorem}
  Let  $V$  be a toric variety defined by a proper simplicial
  fan $\D$ and $X$  an irreducible divisor on $V$  such that 
  $K_V+X$  is 
  terminal.
  Then there exists a sequence of birational toric maps:
$$V=V_1 \stackrel{\varphi_1}{-\to} V_2 \stackrel{\varphi_2}{-\to}
  \cdots \stackrel{\varphi_{r-1}}{-\to}V_r$$
where 

(i) each $\varphi_i$  is either a divisorial contraction or 
  a flip, in particular $V_i$  is defined by 
  a proper simplicial fan;

(ii) for the proper transform $X_i$  of  $X$ on $V_i$ $(i=1,\ldots
  , r)$,
  $K_{V_i}+X_i$  is terminal;

(iii) either that $K_{V_r}+X_r$  is nef or 
  that there exists an extremal ray $R$ on $V_r$ such that 
  $(K_{V_r}+X_r)R<0$ and the elementary contraction $\varphi_R:
  V_r\to Z$  is a fibration to a lower dimensional variety $Z$.
\end{thm}

\begin{pf}
  If $K_V+X$  is nef,
  then the statement is obvious.
  If $K_V+X$  is not nef,
  then there is an extremal ray $R$ such that $(K_V+X)R<0$.
  Take the elementary contraction $\varphi_R:V\to V'$.
  If $\dim V'< \dim V$, then the statement holds.
  So assume that $\varphi _R$  is birational.
  If $\varphi_R$ is divisorial, then define 
  $\varphi_1:=\varphi_R:V\to V'=:V_2$.
  If $\varphi_R$  is not divisorial,
  then let $\varphi_1:V-\to V_2$ be the flip.
  Then in both cases, $K_{V_2}+X_2$  is terminal 
  by Lemma \ref{elementary}.
  Now if $K_{V_2}+X_2$  is nef,
  then the proof is completed.
  If it is not nef,
  make the same procedure as above.
  By the successive procedure,
  one obtains a sequence of divisorial contractions and flips:
$$V=V_1 \stackrel{\varphi_1}{-\to} V_2 \stackrel{\varphi_2}{-\to}
  \cdots \stackrel{\varphi_{r-1}}{-\to}V_r \cdots .$$
  It is sufficient to prove that the sequence terminates at 
  finite stage.
  Let us assume that there exists such a  sequence of infinite 
  length.
  Since the divisorial contraction makes the Picard number 
  strictly less,
  the number of divisorial contractions in the sequence is finite.
  So we may assume that there is $m_0\in \bN$ such that 
  $\varphi_m$'s are all flips for $m\geq m_0$.
  By (iv) of \ref{reid} the set of one dimensional cones of the
  fan defining $V_m$ ($m\geq m_0$) are common.
  As the number of such fans is finite,
  there are numbers $m<m'$ such that 
  $\varphi_{m'-1}\circ \cdots \circ \varphi_m: V_m -\to V_{m'}$
  is identity.
  For each flip $\varphi_j$  $(j=m,\ldots , m'-1)$, 
  take the dominating variety $V_j'$  as in (iv) of \ref{reid}:
 $$\matrix
  &              &    V'_j  &                & \\
  &\psi_j \swarrow &          &\searrow \psi'_j & \\
V_j &              &          &                &V_{j+1}
\endmatrix.$$
  Let $D_j$  be the exceptional divisor of 
  $\psi_j$ and $\psi_{j+1}$.
  Then take a proper toric variety $\tilde V=T_N(\tilde \D)$ 
  which dominates all $V'_j$, $j=m,\ldots , m'-1$ and
  on which the proper transform $\tilde X$  of  $X_j$'s   is 
  $\tilde \D$-regular. 
  This is possible, because $K_{V_j}+X_j$'s are terminal.
  Here one should note that the set of exceptional divisors
  on $\tilde V$  for all morphisms $\tilde V\to V_j$ 
  $(j=m,\ldots ,m'-1)$ are common.
  For every $j=m,\ldots , m'-1$,
  the discrepancy $\alpha$ of $K_{V_j}+X_j$ at $D_j$  is 
  less than the discrepancy $\alpha'$  of 
  $K_{V_{j+1}}+X_{j+1}$ at $D_j$ by \ref{elementary}.
  By this fact, for every exceptional divisor $E$  on $\tilde V$,
  the discrepancy $\alpha_E$ of $K_{V_j}+X_j$ at $E$ 
  and  the discrepancy $\alpha'_E$ of $K_{V_{j+1}}+X_{j+1}$ at 
  $E$  satisfy $\alpha_E\leq \alpha'_E$
   and for 
  at least one exceptional divisor $E$, $\alpha_E<\alpha'_E$.
  Therefore comparing   $K_{V_m}+X_m$
  and  $K_{V_{m'}}+X_{m'}$,
  there exists an exceptional divisor on $\tilde V$
  at which the discrepancy of $K_{V_m}+X_m$ is less than that of 
  $K_{V_{m'}}+X_{m'}$, which is the contradiction to that 
  $V_m\to V_{m'}$  is the identity.
\end{pf}

  To apply the theorem above to the minimal model problem for a 
  toric divisor, 
  one needs the following lemma.

\begin{lem}
\label{basic lemma}
(Lemma 2.7, \cite{weight})
 Let  $Y\subset Z$  be an irreducible  Weil divisor on a variety   
   $Z$. Assume that  $Z$ admits at worst $\bQ$-factorial log-terminal
singularities.
  Let  $\varphi:\tilde Y \to Y$  be a resolution of singularities on  $Y$.
  Assume $K_{\tilde Y}=\varphi^*((K_Z+Y)|_{Y})+\sum_im_i E_i$
  with  $m_i>-1$  for all $i$,
  where $E_i$'s are the exceptional divisors of  $\varphi$.

   Then  $Y$  is normal,
  and $Y$  has at worst log-terminal singularities.
  
  In particular,
if $m_i> 0$ for all $i$,
then  $Y$  has at worst terminal singularities.
\end{lem}

\begin{cor}
\label{MMT}
  Let  $V$  be a toric variety defined by a proper fan $\D$
  and $X$ a $\D$-regular divisor on $V$.
  If $\kappa(X)\geq 0$, then $X$  has a minimal model with the
  abundance.
  If $\kappa(X)=-\infty$, then $X$  is birationally equivalent
  to a proper variety $Y$ with at worst terminal singularities
  and a fibration $\varphi:Y\to Z$ to a lower dimensional 
  variety $Z$ with $-K_{Y}$ relatively ample.
\end{cor}

\begin{pf}
  Let $V_1$  be the toric variety defined by a non-singular
  subdivision  $\D_1$ of  $\D$ and $X_1$ be the proper transform
  of $X$  on  $V_1$.
  Then  $X_1$  is $\D_1$-regular
  and therefore $K_{V_1}+X_1$  is terminal by \ref{terminal}.
  Then  one obtains a sequence: 
 $$V_1 \stackrel{\varphi_1}{-\to} V_2 \stackrel{\varphi_2}{-\to}
  \cdots \stackrel{\varphi_{r-1}}{-\to}V_r$$
  as in  Theorem \ref{theorem}.
  One can prove that for each $j=1, \ldots , r$, $X_j$ has at worst
  terminal singularities.
  In fact, take a morphism $\varphi:\tilde V \to V_j$
  corresponding to a non-singular subdivison $\tilde \D$ of the 
  fan $\D_j$ of $V_j$ such that the proper transform $\tilde
  X$  of $X$ is $\tilde \D$-regular.
  Then, as $K_{V_j}+X_j$ is terminal, it follows that
$$(K_{\tilde V}+\tilde X)|_{\tilde X}=
  \varphi^*((K_{V_j}+X_j)|_{X_j})+\sum_ia_iE_i|_{\tilde X}
  \ \ \ \ (a_i>0\ \text{for\ all}\ i).$$
  Here the left hand side is the canonical divisor $K_{\tilde X}$
  of a non-singular variety $\tilde X$.
  Therefore by Lemma \ref{basic lemma}, one sees that $X_j$  has
  at worst terminal singularities.
  By (iii)  of \ref{theorem} there are two cases for $V_r$.
\begin{case}
  $K_{V_r}+X_r$ is nef.

  Then the linear system $|m(K_{V_r}+X_r)|$  is basepoint
  free for some $m\in \bN$.
  This is proved by a slight modification of the proof of
  Toric Nakai Criterion
  (2.18, \cite{Oda88}).
  Therefore $|mK_{X_r}|$  is basepoint free,
  which implies that $X_r$  is a minimal model with the abundance.
  In this case,  $\kappa(X)=\kappa(X_r)\geq 0$.
\end{case}
\begin{case}
  There exists an extremal ray $R$ on $V_r$ such that 
  $(K_{V_r}+X_r)R<0$ and the elementary contraction $\varphi_R:
  V_r\to Z$  is a fibration to a lower dimensional variety $Z$.

  Under this situation, first consider the case:

{\sl Subcase.}  
  $\dim X_r> \dim \varphi_R(X_r)$.

  Let $ F$  be a fiber  of 
   $\varphi_R$.
  Then by (ii) of \ref{reid}, $F$  is a weighted projective
  space and $(K_{V_r}+X_r)C<0$ for every curve $C$  in $F$,
  which implies that $-(K_{V_r}+X_r)$ is relatively ample
  over $Z$.
  Hence $-K_{X_r}$ is relatively ample over $\varphi_R(X_r)$.
  This yields that $\kappa(X)=\kappa(X_r)=-\infty$,
  and $\varphi_R|_{X_r}:X_r\to \varphi_R(X_r)$  is a desired
  fibration.

{\sl Subcase.} $\dim X_r=\dim \varphi_R(X_r)$.

  In this case $\dim Z=\dim V_r-1$ and every fiber $\ell$ 
  of $\varphi_R:
  V_r \to Z$  is $\bP^1$ by (ii) of \ref{reid}.
  Therefore  $K_{V_r}\ell=-2$.
  On the other hand, because $\varphi|_{X_r}$ is generically 
  finite, $X_r\ell>0$.
  Here, since $V_r$ has at worst terminal singularities by 
  \ref{terminal sing}, the singular locus has codimension greater
  than 2 and therefore the divisor $X_r$ is a Cartier divisor
  along a general fiber  $\ell$,
  which yields that $X_r\ell$  is an integer.
  By $(K_{V_r}+X_r)\ell<0$, it follows $X_r\ell=1$.
  It implies that $\varphi_R|_{X_r}:X_r\to Z$  is 
  a birational morphism,
  therefore $X_r$  is rational.
  So $X$ and $X_r$  are birationally equivalent to 
  $\bP^n$ which has ample anti-canonical divisor 
  and of course $\kappa(X)=-\infty$.
\end{case}
\end{pf}

\begin{cor}
  Let the ground field $k$ be of characteristic zero.
  Let $V$ be  a proper toric variety,  $|L|$ a linear system
   without a basepoint and $X$  a general member of $|L|$.
  Then the statements of Corollary \ref{MMT} hold for $X$.
\end{cor}
\begin{pf}
By the Bertini's Theorem, $X$  is $\D$-regular.
\end{pf} 

\begin{cor}
\label{kappa}
  Let $V$  be a toric variety defined by a proper fan $\D$
  and $X$ a $\D$-regular divisor on $V$.
  Assume $\kappa(X)\geq 0$.
  Then there exists a non-singular
  subdivision $\tilde \D$ of $\D$
  such that  $\tilde V=T_N(\tilde \D)$  and 
  the proper transform $\tilde X$ of $X$ on $\tilde V$
  satisfy the following: 
 $$\kappa(\tilde V, K_{\tilde V}+\tilde X)\geq 0.$$
\end{cor}

\begin{pf}
  Use the notation of the proof of  \ref{MMT}.
  Take a nonsingular subdivision $\tilde \D$ of 
  both $\D$ and  $\D_r$ which is the fan of $V_r$.
  Then the proper transform $\tilde X$  of $X$ on 
  $\tilde V=T_N(\tilde \D)$
  is $\tilde \D$-regular.
  Since $K_{V_r}+X_r$  is terminal and 
  $|m(K_{V_r}+X_r)|$ is basepoint free for some $m\in \bN$,
$$0\neq \Gamma(V_r, m(K_{V_r}+X_r))\subset 
   \Gamma(\tilde V, m(K_{\tilde V}+\tilde X)).$$
\end{pf}

\section{\bf Divisors and Polytopes}

\begin{say}
  Here we summerize the basic notion of an invariant divisor of a toric
  variety and the corresponding polytope which will be used in the 
  next section.
  In this paper,
  a polytope in an $\bR$-vector space means the intersection of finite
  number of 
  half-spaces  $\{\bm | f_i(\bm)\geq a_i\}$ for linear functions $f_i$.
\end{say}

\begin{say}
  Let  $M$  be the free abelian group  ${\bZ}^n$ $(n\geq 3)$
  and  $N$  be the dual $Hom_{\bZ}(M, {\bZ})$.
  We denote  $M\otimes _{\bZ}{\bR}$  and $N\otimes_{\bZ}{\bR}$  by
  ${M_{\bR}}$  and  $\NR$, respectively.
  Define  $M_{\bQ}$  and  $N_{\bQ}$  in the same way.
  Then one has the canonical pairing  $(\ \ ,\ \ ):N\times M \to \bZ$,
  which can be canonically extended to 
  $   (\ \ ,\ \ ):\NR\times M_{\bR} \to \bR$.
  For a  fan ${\D}$ in ${\NR}$,
  we construct the toric variety  $T_N({\D})$.
  The fan $\D$  is always assumed to be proper,
i.e. the support  $|\D|=\NR$.
  Denote by $\D(k)$ the set of $k$-dimensional cones in $\D$.
  Denote by  ${\D}[1]$
  the set of primitive vectors ${\bq}=(q_1,\ldots , q_r)\in {N}$
  whose rays  ${\bR}_{\geq 0}{\bq}$  belong to ${\D(1)}$.
  For ${\bq}\in {\D}[1]$,
  denote by $D_{\bq}$ the
  corresponding divisor which is denoted by $\overline{orb\ 
  {\bR}_{\geq 0}{\bq}}$  in  \cite{Oda75}.
  Denote by $U_{\sigma}$  the invariant affine open subset which contains
  ${orb\ 
  {\sigma}}$ as the unique closed orbit.
\end{say}
\begin{defn}
  For $\bp \in \NR$ and  a subset $K\subset \MR$,
  define
  $$\bp(K):=\inf_{\bm\in K}(\bp, \bm)$$
\end{defn}

\begin{defn}
\label{support function}
  Let $\D$  be a proper  fan in  $\NR$.
  A continuous function  $h:\NR \to \bR$  is called
  a $\D$-support function, if 

  (1) $h|_{\sigma}$ is $\bR$-linear for every cone  $\sigma \in \D$
  and

  (2)  $h$  is  $\bQ$-valued on  $N_{\bQ}$.

  A $\D$-support function  $h$  is called integral if 

(2')  $h$  is $\bZ$-valued on $N$.
\end{defn}
\vskip.5truecm
\begin{prop}
\label{correspondence}
  For a $\D$-support function  $h$,
  define  $D_h=-\sum_{\bp\in \D[1]}h(\bp)D_{\bp}$.
  Then the correspondence  $h \mapsto D_h$ gives a bijective map:

  \{$\D$-support functions\} $\simeq$ \{invariant  $\bQ$-Cartier divisors 
  on  $\tnd$\}.
  Here  $D_h$  is a Cartier divisor,
  if and only if  $h$  is integral.
\end{prop}
\vskip.5truecm
\begin{defn}
  For a $\D$-support function  $h$,  define
$$\boxh:=\{m\in M_{\bR}|(\bp,\bm)\geq h(\bp), \forall \bp\in \NR\},$$
  and call it the polytope associated with $h$ or with $D_h$.
  Actually it is a polytope by \ref{boundary} and compact since 
  the fan $\D$  is proper.
\end{defn}
\vskip.5truecm
\begin{prop}
\label{basepoint free}
  (see \cite{Oda88})
  For an integral $\D$-support function $h$, 
  the following are equivalent:

  (i) the linear system $|D_h|$  is basepoint free;

  (ii)  $h$  is upper convex; i.e. for arbitrary  $\bn, \bn'\in \NR$,
   $h(\bn)+h(\bn')\leq h(\bn+\bn')$;

   (iii) $\boxh=$ the convex hull of $\{h_{\sigma}|\sigma \in \D(n)\} $,
   where $h_{\sigma}$  is a point of $M$ which gives the linear
  function
  $h|_{\sigma}$  for $\sigma \in \D(n)$.
\end{prop}
\vskip.5truecm
\begin{prop}
\label{ample}
  (see \cite{Oda88})
  For a $\D$-support function $h$, 
  the following are equivalent:

  (i) the $\bQ$-Cartier divisor $D_h$ is ample;

  (ii)  $h$  is strictly upper convex; i.e. $h$ is upper convex and 
   $h(\bn)+h(\bn')< h(\bn+\bn')$, if there is no  cone $\sigma$ such that
   $\bn, \bn'\in \sigma$;

   (iii) $\boxh$ is of dimension  $n$ and the correspondence
  $\sigma \mapsto h_{\sigma}$ gives the bijective map 
  $\D(n) \simeq $ \{The vertices of $\boxh$\},
     where $h_{\sigma}$  is a point of $M_{\bQ}$ which gives the linear 
  function
  $h|_{\sigma}$  for $\sigma \in \D(n)$.
\end{prop} 
\vskip.5truecm
  Now we show simple lemmas which are used in the next section.
\begin{lem}
\label{inf}
  Let $h$  be a $\D$-support function.
  If  $h_{\sigma}\in \boxh$  for every  $\sigma\in \D(n)$,
  then $h(\bp)=\bp(\boxh)$  for every  $\bp\in \NR$, and
  the polytope  $\boxh$  is the convex hull of the set $\{h_{\sigma}\}$.
\end{lem}
\begin{pf}
  By the definition of $\boxh$,
  $h(\bp)\leq(\bp,\bm)$ for all $\bm\in \boxh$.
  Therefore  $h(\bp)\leq \bp(\boxh)$.
  Let  $\sigma$  be the cone in $\D(n)$  such that  $\bp\in \sigma$,
  then $h(\bp)=(\bp, h_{\sigma})\geq \bp (\boxh)$, since $h_{\sigma}\in \boxh$. 
  For the second assertion,
  assume a vertex $\bm\in \boxh$ does not belong to the convex hull of 
  $\{h_{\sigma}\}$.
  Then there exists $\bp\in \NR$ such that  $(\bp, \bm)<(\bp,h_{\sigma})$
  for every  $\sigma \in \D(n)$,
  where the left hand side is greater than or equal to $h(\bp)$
  by the definition of  $\boxh$.
  This is a contradiction,
  because for $\sigma\in \D(n) $ such that  $\bp \in \sigma$,
  $h(\bp)=(\bp, h_{\sigma})$.   
\end{pf}
\vskip.5truecm
\begin{lem}
\label{boundary}
  Denote an invariant divisor  $D_h=\sum_{\bp\in \D[1]}m_{\bp}D_{\bp}$.
  Then $$\boxh=\bigcap_{\bp\in \D[1]}\{\bm\in \MR|
  (\bp,\bm)\geq -m_{\bp}\}.$$
\end{lem}
\begin{pf}
  By \ref{correspondence},
  $m_{\bp}=-h(\bp)$,
  then the inclusion  $\boxh\subset \bigcap_{\bp\in \D[1]}\{\bm\in \MR|
  (\bp,\bm)\geq -m_{\bp}\} $  is obvious.
  Take an element $\bm$ from the right hand side.
  For an arbitrary $\bp \in \NR$,
  take  $\sigma \in \D(n)$ such that  $\bp \in \sigma$.
  Let $\sigma$  be spanned by  
  $\bp_1,\bp_2,\ldots ,\bp_s$ $(\bp_i\in \D[1])$,
 then $\bp=\sum a_i\bp_i$ with $a_i\geq 0$.
  One obtains that   
  $(\bp,\bm)= \sum a_i(\bp_i,\bm)\geq \sum a_ih(\bp_i)
  =\sum a_i(\bp_i,h_{\sigma})=h(\bp)$,
  which shows that $\bm$ belongs to $\boxh$.
\end{pf}
\vskip.5truecm
\begin{defn}
\label{contribution}
  Let $\Box$  be a polytope in $\MR$ defined by  
  $\bigcap_{i=1}^rH_i$,  where
  $H_i=\{\bm\in \MR| (\bp_i,\bm)\geq a_i\}$. 
  We say that  $H_i$ contributes to $\Box$,  
  if $\Box \cap \{\bm\in \MR| (\bp_i,\bm)= a_i\}\neq \phi$.
  And we say that $H_i$ contributes properly to $\Box$,
  if $\bigcap_{j\neq i}H_j\neq \Box$.
\end{defn}

\begin{defn}
\label{dual fan}
  Let  $\Box$  be an $n$-dimensional compact polytope in  $\MR$.
  Define the dual fan $\Gamma_{\Box}$ of $\Box$ as follows:
  $\Gamma_{\Box}=\{\gamma^*\}$, where $\gamma$  is a face of  $\Box$
  and  $\gamma^*:=\{\bn\in\NR|$ the function $\bn|_{\Box}$ attains the
minimal value
  at all points of $\gamma$\}.
  Then  $\Gamma_{\Box}$  turns out to be a proper fan.
\end{defn}
\vskip.5truecm
\begin{say}
\label{projective}
  If  $\D$  is the dual fan of the polytope  $\boxh$ corresponding to
  a $\D$-support function $h$,
  then by \ref{ample} $D_h$  is ample,
  therefore the variety $\tnd$  turns out to be a projective
  variety.
\end{say}

\section{\bf The construction of a minimal model}

\begin{say}
  In this section we concretely construct 
  a projective minimal model
  with the abundance  for a $\D$-regular toric divisor $X$ with
  $\kappa(X)\geq 0$
  by means of a polytope of the adjoint divisor.
  Let $V$ be a toric variety defined by a proper fan $\D$
  and $X$ a $\D$-regular divisor with $\kappa(X)\geq 0$.
  To construct a minimal model of $X$ we may assume that
  $V$  is non-singular and $\kappa(V, K_V+X)\geq 0$,
  by Corollary \ref{kappa}.
\end{say}

\begin{say} 
\label{construction}
  {\bf The construction}
  Let  $h$  be a $\D$-support function such that $ K_{\tnd}+X \sim D_h$.
  Then, by $\kappa (\tnd , K_{\tnd}+X) \geq 0$, 
  it follows that  $\boxh \neq \phi$.
  Let $\D[1]=\{\bp_1,\ldots ,\bp_s\}$ and $H_i=\{\bm \in \MR|(\bp_i,\bm)
  \geq h(\bp_i)\}$.
  Then by \ref{boundary}
  $\boxh=\bigcap_{i=1}^sH_i$.
  Assume that $H_1,\ldots ,H_r$ $(r\leq s)$ are all that contribute 
  to $\boxh$.
  For   $\epsilon_i>0$  $i=1,\ldots , r$,
  define 
  $H_{i, \epsilon _i}
  :=\{\bm\in \MR| (\bp_i,\bm)\geq h(\bp_i)-\epsilon_i\}$,
  $\partial H_{i, \epsilon _i}
  :=\{\bm\in \MR| (\bp_i,\bm)= h(\bp_i)-\epsilon_i\}$
  and $\Box({\epsilon}):=\bigcap_{i=1}^rH_{i,\epsilon_i}$,
  where  $\epsilon = (\epsilon_1,\ldots , \epsilon_r)$.
  Here one should note that the polytope $\boxh$ may not be 
  of the maximal dimension.
  By "puffing up" this, one get a polytope $\Box({\epsilon})$
  of the maximal dimension.
  The subset  $Z=\{\epsilon \in \bR_{>0}^r| \bigcup
  \partial H_{i,\epsilon_i}$
  is not of normal crossings\}    is Zariski closed
   and the complement $\bR_{>0}^r\setminus Z$  is divided 
  into finite number of chambers.
  Take a chamber  $W$  such that:

  (\ref{construction}.1) $0\in \overline{W}$;

  (\ref{construction}.2) every $H_{i,\epsilon_i}$  $ (i=1,\ldots, r)$ 
  contibutes properly to
  $\Box({\epsilon})$ for  $\epsilon\in W$. 

 Then the dual fan $\Sigma$ of  $\Box({\epsilon})$  is common for 
  every  $\epsilon \in W$ and it is simplicial,
  because $\bigcup\partial H_{i,\epsilon_i}$  is of normal crossings.
  Let  $\xs$  be the proper transform of $X$  in  $\tns$.
  we claim that $\xs$ is a minimal model of  $X$ with the abundance.
  One can see that $\tns$ is projective,
  because an invariant $\bQ$-Cartier divisor  
  $\sum_{\bp_i\in\Sigma[1]}(h(\bp_i)-\epsilon_i)D_{\bp_i}$
  with all $\epsilon_i$ rational and $\epsilon \in W$
  is  ample
  since  $\Sigma$ is the dual fan of the corresponding polytope
  to this divisor (\ref{projective}).
  Hence the projectivity of $\xs$ follows automatically.
\end{say}
\vskip.5truecm
\begin{say}
  Now we are going to prove that $\xs$ satisfies desired conditions
  for a minimal model.
  First note that $\Sigma[1]=\{\bp_1,\ldots , \bp_r\}$,
  by \ref{construction}.2.
  Next note that every $\bQ$-Weil divisor on $\tns$  is a $\bQ$-Cartier
  divisor,
  because $\Sigma$ is simplicial and therefore 
  $\tns$  has quotient singularities.
\end{say}
\begin{claim}
\label{h and k} 
  The divisor  $K_{\tns}+\xs$ is linearly equivalent to an invariant
  divisor  $-\sum_{i=1}^rh(\bp_i)D_{\bp_i}$.
  Let  $k$  be the $\Sigma$-support function corresponding to 
  this divisor,
  then $h(\bp_i)=k(\bp_i)$  for $i=1,\ldots , r$ and $\boxh=\boxk$.
\end{claim}
\begin{pf}
  The first assertion follows from that the divisor
  $K_{\tns}+\xs$  is the proper transform of $K_{\tnd}+X \sim 
  -\sum_{i=1}^sh(\bp_i)D_{\bp_i}$.
  The second assertion is obvious and the last assertion
  follows from \ref{boundary}
  and the fact that $H_1,\ldots , H_r$ are all that contribute to
  $\boxh$.
\end{pf}
\begin{claim}
\label{k in box}
  For all $\sigma \in \Sigma(n)$, 
   it follows that  $k_{\sigma}\in \boxk$.
\end{claim}
\begin{pf}
  Let  $\{\epm\}_m$  be  series of rational points in $W$
  which converge to $0$.
  Let  $\km$  be the $\Sigma$-support function corresponding to
  a $\bQ$-Cartier divisor 
  $\sum_{i=1}^r(-h(\bp_i)+\epm _i)D_{\bp_i}$.
  Then by \ref{boundary} 
  it follows that $\Box_{\km}=\Box(\epm)$, 
  and therefore by \ref{projective}
  the divisor is ample.
  Replacing $\{\epm\}_m$  by suitable subsequence,
  one can assume there exists $\lim_{m\to \infty}\km_{\sigma}$
  for every $\sigma\in \Sigma (n)$.
  Indeed,  replacing by suitable subsequence,
  one may assume that $\epm_i\geq {\epsilon^{(m+1)}}_{i}$ for every $i$,
  then $\Box_{\epm}\supset \Box_{\epsilon^{m+1}}\supset\cdots$;
  therefore for every  $\sigma\in \Sigma(n)$  and $m$
  it follows that  $\km_{\sigma}\in
  \Box_{\epsilon^{(1)}}$ which is compact;
  so $\{\km_{\sigma}\}$ have  an accumulating point. 
  Let $k'_{\sigma}:=\lim_{m\to \infty}\km_{\sigma}$,
  then $k'_{\sigma}\in \boxk$,
  because the ampleness of $D_{\km}$ yields
  $\km_{\sigma}\in \Box(\epm)$.
  The collection $\{k'_{\sigma}\}_{\sigma\in\Sigma(n)}$ defines a
   function $k'$ on $\NR$.
  In fact, for every $m$,  $\km_{\sigma}=\km_{\tau}$ as a function on
  $\sigma\cap \tau$,
  which yields that $k'_{\sigma}=k'_{\tau}$  as a function on 
  $\sigma\cap \tau$.
  Now one obtains that $k'=k$.
  This is proved as follows:
  for every $\bp_i\in \Sigma[1]$ take  $\sigma\in \Sigma(n)$
  such that $\bp_i \in \sigma$;
  $k'(\bp_i)=(\bp_i,k'_{\sigma})=\lim_{m\to \infty}(\bp_i,\km_{\sigma})
  = \lim_{m\to \infty}(h(\bp_i)-\epm_i)=h(\bp_i)=k(\bp_i)$,
  since  $\km_{\sigma}$  is on the hyperplane 
  $(\bp_i,\bm)=h(\bp_i)-\epm_i$.
  Hence it follows that $k=k'$ and therefore $k_\sigma =k'_\sigma $  
  for every $\sigma\in \Sigma (n)$, 
  which shows that $k_{\sigma}\in \boxk$.
\end{pf}

  Now by \ref{inf} and \ref{basepoint free} 
  the linear system $|mD_k|=|m(K_{\tns}+\xs)|$ has no basepoint
  for such $m$  that  $mD_k$  is a Cartier divisor.
\begin{say}
  Let  $\tsig$  be a non-singular subdivision of $\Sigma$  and $\D$.  
  Let 
$$\matrix
  & \psi \nearrow & \tnd\\
T_N(\tsig)& & \\
  &  \varphi \searrow & \tns
\endmatrix$$
be the corresponding morphisms and  $X(\tsig)$ the proper transform
  of  $X$  in  $T_N(\tsig)$.
  Since  $X(\tsig)$ is $\tsig$-regular by \cite{Batyrev},
  it is non-singular and $\varphi |_{X(\tsig)}$ 
  is birational.
\end{say}
\begin{claim}
\label{positive}
  It follows that
  $$K_{T_N(\tsig)}+X(\tsig)=\varphi^*(K_{\tns}+X(\Sigma))+
  \sum_{\bp\in \tsig[1]\setminus \Sigma[1]}m_{\bp}D_{\bp},$$
  where $m_{\bp}>0 $  for $\bp$  such that  
  $D_{\bp}\cap X(\tsig)\neq\phi$.
\end{claim}
\begin{pf}
  Denote
  $$K_{T_N(\tsig)}+X(\tsig)=\psi^*(K_{\tnd}+X)+\sum_{\bp\in\tsig[1]
  \setminus\D[1]}\alpha_{\bp}D_{\bp},$$
  then $\alpha_\bp>0$  for $\bp$  such that  $D_{\bp}\cap X(\tsig)\neq\phi$,
  since $X$  is non-singular.
  Putting $\alpha_\bp=0$ for $\bp\in\D[1]$,
  one obtains that 
  $K_{T_N(\tsig)}+X(\tsig)\sim \sum_{\bp \in \tsig[1]}
  (-h(\bp)+\alpha_{\bp})D_{\bp}$,
  as  $K_{\tnd}+X\sim D_h$.
  On the other hand,
  $$K_{T_N(\tsig)}+X(\tsig)=\varphi^*(K_\tns+X(\Sigma))+
  \sum_{\bp\in \tsig [1]\setminus \Sigma [1]}m_\bp D_\bp.$$
Putting $m_\bp=0$ for $\bp\in\Sigma[1]$,
  one obtains that 
  $K_{T_N(\tsig)}+X(\tsig)\sim \sum_{\bp \in \tsig[1]}
  (-k(\bp)+m_{\bp})D_{\bp}$,
  as  $K_{\tns}+X(\Sigma)\sim D_k$.

  Therefore 
$\sum_{\bp \in \tsig[1]}
  (-h(\bp)+\alpha_{\bp})D_{\bp}\sim
  \sum_{\bp \in \tsig[1]}
  (-k(\bp)+m_{\bp})D_{\bp}$.
  As  $h(\bp)=k(\bp)$ and $\alpha_\bp=m_\bp=0$ for  $\bp\in \Sigma[1]$,
  one obtains that
  $$\sum_{\bp\in\tsig[1]\setminus\Sigma[1]}((-h(\bp)+\alpha_\bp)-
  (-k(\bp)+m_\bp))D_\bp \sim 0.$$
  Here $D_\bp$  $(\bp\in\tsig[1]\setminus\Sigma[1])$ 
  are all exceptional for $\varphi$.
  Then the divisor above is not only linearly equivalent to 0
  but also equal to 0.
  Therefore $(-h(\bp)+\alpha_\bp)-
  (-k(\bp)+m_\bp)=0$ for every $\bp\in \tsig[1]\setminus\Sigma[1]$,
  where $k(\bp)=\bp (\boxh)$ by $k_\sigma \in \boxk=\boxh$ and 
  \ref{inf}.
  Now consider the divisor $D_\bp$  such that  
  $D_\bp\cap X(\tsig)\neq 0$.
  For $\bp\in \tsig[1] \setminus \D[1]$,
  $m_\bp=\bp(\boxh)-h(\bp)+\alpha_\bp\geq \alpha_\bp > 0$.
  For $\bp\in \D[1]\setminus\Sigma[1]$,
  it follows that $m_\bp=\bp(\boxh)-h(\bp)>0$,
  because  $\{\bm|(\bp,\bm)\geq h(\bp)\}$  does not contribute to
  $\boxh$ by the definition of $\Sigma$ (c.f \ref{construction}).
  This completes the proof.
\end{pf}
\vskip.5truecm
\begin{say}
  Since  $\tns$ has at worst quotient singularities,
  one can apply Lemma\ref{basic lemma} to our situation
  and obtain that $\xs$  is normal and has at worst terminal
  singularities. And the linear system of  
  $mK_\xs=m(K_{\tns}+\xs)|_{\xs}$ $(m\gg 0)$ has no basepoint,
  because  $|m(K_{\tns}+\xs)|$ is basepoint free
  as is noted after the proof of \ref{k in box}.
  This completes the proof of that $\xs$ is a projective minimal
  model with the abundance.
\end{say}

\begin{say}
  Pursuing elementary contractions and flips is like groping
  for a minimal model in the dark.
  The reason why the discussion of this section goes well
  without contractions nor flips is because in toric geometry
  every exceptional divisor is visible as a vector in 
  the space $N$.
  Then one can prepare so that every discrepancy of adjoint divisor
  is positive (cf. \ref{positive}),
  which makes the singularities terminal.
  In the discussion, one puffed up the polytope of the adjoint
  divisor and took its dual fan $\Sigma$.
  This implies that in $\tns$ the adjoint divisor is the limit
  of a sequence of ample divisors (cf. \ref{k in box}),
  which makes the adjoint divisor nef;
  or equivalently semi-ample.
\end{say}

\vskip 1truecm
\section{\bf Examples}

  In this section the base field $k$ is always assumed to be of 
  characteristic zero.
  Let  $M$  be  $\bZ^3$ and $N$  be its dual.

\begin{exmp}
\label{exmp1}
  Let $\bp_i$ $(i=1,\ldots , 6)$  and  $\bq_j$  $(j=1, \ldots , 8)$
  be points in $N$  as follows:
  $\bp_1=(1,0,0),$ $\bp_2=(-1,0,0),$  $\bp_3=(0,1,0),$ $\bp_4=(0,-1,0),$
  $\bp_5=(0,0,1),$ $\bp_6=(0,0,-1),$
  $\bq_1=(1,1,1),$ $\bq_2=(-1,-1,-1),$ $\bq_3=(1,1,-1),$ $\bq_4=(-1,-1,1),$
  $\bq_5=(1,-1,1),$ $\bq_6=(-1,1,-1),$ $\bq_7=(-1,1,1),$ $\bq_8=(1,-1,-1).$
  Let them generate one-dimensional cones  $\bR_{\geq 0}\bp_i$, 
  $\bR_{\geq 0}\bq_j$ and construct a fan  $\D$ with these cones
  as in Figure 1.
  Here note that Figure 1 is the picture of the fan which is cut 
  by a hypersphere with the center the origin and unfolded onto
  the plane.
  This fan is the dual fan of the polytope of Figure 2 and it is easy to
  check that it is non-singular.
  Let  $X$  be a general member of a base-point-free linear system
  $|\sum_{i=1}^6D_{\bp_i}
  +2\sum_{j=1}^8D_{\bq_j}|$.
  Let $h$  be the $\D$-support function such that  $K_{\tnd}+X\sim D_h$.
  Then the polytope  $\boxh$  is one point 
  $\bigcap_{i=1}^6\{\bm| (\bp_i,\bm)\geq 0\}$and the half spaces 
  contributing
  to this polytope are $\{\bm| (\bp_i,\bm)\geq 0\}$, $i=1,\ldots ,6$,
  because $K_{\tnd}+X\sim \sum D_{\bq_j}$.
  Therefore, for a sufficiently small general $\epsilon$,
  the polytope 
  $\Box(\epsilon)=\bigcap_i\{\bm\in \MR|(\bp_i,\bm)\geq -\epsilon_i\}$
  is a hexahedron 
   whose picture is as Figure 3.
  The dual fan $\Sigma$ of $\Box(\epsilon)$ (Figure 4) gives a 
  minimal model $\xs$ of $X$.
  Since  $K_{\tns}+\xs \sim 0$,
  it follows that $\kappa (X)= 0$.
\end{exmp}

\begin{exmp}
  Let $\bp_i$  and $\bq_j$  be as in \ref{exmp1}
  and  $\D$  be the fan with the cones generated by these vectors
  as Figure 5.
  This fan is the dual fan of the polytope of Figure 6 and it is 
  easy to check that it is non-singular.
  Let  $X$  be a general member of a base-point-free linear system
  $|2D_{\bp_1}+2D_{\bp_2}+\sum_{i=3}^6D_{\bp_i}+3\sum_{j=1}^8D_{\bq_j}|$.   
  Let $h$  be the $\D$-support function such that  $K_{\tnd}+X\sim D_h$.
  Then the polytope  $\boxh$  is a segment
  $\bigcap_{i=1}^2\{\bm|(\bp_i,\bm)\geq -1\}\cap \bigcap_{i=3}^6
  \{\bm | (\bp_i,\bm)\geq 0\}$ 
  and the half spaces contributing 
  to this polytope are $\{\bm|(\bp_i,\bm)\geq -1\}$ $(i=1,2)$
  and $\{\bm | (\bp_i,\bm)\geq 0\}$ $(i=3,\ldots ,6)$,
  because $K_{\tnd}+X\sim D_{\bp_1}+D_{\bp_2}+ 2\sum D_{\bq_j}$.
  Therefore, for a sufficiently small general $\epsilon$,
  the polytope 
  $\Box(\epsilon)=
   (\bigcap_{i=1}^2\{\bm\in \MR|(\bp_i,\bm)\geq -1-\epsilon_i\})
  \cap 
  (\bigcap_{i=3}^6\{\bm\in \MR|(\bp_i,\bm)\geq -\epsilon_i\})$
  is a hexahedron 
   whose picture is as Figure 7.
  The dual fan $\Sigma$ of $\Box(\epsilon)$ (Figure 4) gives a 
  minimal model $\xs$ of $X$.
  Since $\boxh$ is of one dimension,
  $\dim \Gamma(\tns , m(K_{\tns}+\xs ))$ grows in order 1,
  and therefore $\dim \Phi _{|m(K_{\tns}+\xs )|}(\tns)=1$.
  This shows that $\dim \Phi_{|mK_{\xs}|}(\xs)\leq 1$.
  As the dual fan of the polytope of $\xs \sim 2D_{\bp_1}+2D_{\bp_2}+
  \sum_{i=3}^6D_{\bp_i}$ is $\Sigma$,
  $\xs$  is ample by \ref{projective}.
  Hence  $\xs$ intersects all fibers of  $\Phi _{|m(K_{\tns}+\xs )|}$,
  which shows that $\kappa (X)=1$.
\end{exmp}


\makeatletter \renewcommand{\@biblabel}[1]{\hfill#1.}\makeatother


\vskip 2truecm

\end{document}